\documentclass[11pt]{article}

\usepackage{amsmath, amsfonts, amsthm, amssymb, mathtools}
\usepackage{authblk}
\usepackage{color, bm, graphicx, epstopdf}
\usepackage[font = small, labelfont = bf, labelsep = period]{caption}
\usepackage{subcaption}
\usepackage{xspace}
\usepackage{enumitem}
\usepackage{mathrsfs, bbm}
\usepackage{url}

\usepackage{hyperref}
\usepackage{algorithm}
\usepackage{algpseudocode}
\usepackage{comment}

\usepackage{threeparttable}
\usepackage{multirow}
\usepackage{multicol}
\usepackage{booktabs}
\usepackage{csquotes}
\usepackage{rotating}
\usepackage[round]{natbib}

\graphicspath{{figures/}}

\DeclarePairedDelimiter\ceil{\lceil}{\rceil}

\usepackage{fullpage}[2cm]
\title{A Novel Mixed-Integer Linear Programming Formulation for Continuous-Time Inventory Routing}

\author[1,2]{Akang~Wang}
\author[1]{Xiandong~Li}
\author[3]{Jeffrey~E.~Arbogast}
\author[3]{Zachary~Wilson}
\author[1]{Chrysanthos~E.~Gounaris\thanks{Corresponding author: gounaris@cmu.edu}}
\affil[1]{\small Department of Chemical Engineering and Center for Advanced Process Decision-making, Carnegie Mellon University, Pittsburgh, PA 15213, USA}
\affil[2]{\small Shenzhen Research Institute of Big Data, Shenzhen, Guangdong, China}
\affil[3]{\small Air Liquide, Newark, DE 19702, USA}

\date{}

\begin{document}

\maketitle

\begin{abstract}
	Inventory management, vehicle routing, and delivery scheduling decisions are simultaneously considered in the context of the inventory routing problem. This paper focuses on the continuous-time version of this problem where, unlike its more traditional discrete-time counterpart, the distributor is required to guarantee that inventory levels are maintained within the desired intervals at any moment of the planning horizon. In this work, we develop a compact mixed-integer linear programming formulation to model the continuous-time inventory routing problem. We further discuss means to expedite its solution process, including the adaptation of well-known rounded capacity inequalities to tighten the formulation in the context of a branch-and-cut algorithm. Through extensive computational studies on a suite of $90$ benchmark instances from the literature, we show that our branch-and-cut algorithm outperforms the state-of-the-art approach. We also consider a new set of $63$ instances adapted from a real-life dataset and show our algorithm's practical value in solving instances with up to $20$ customers to guaranteed optimality.

\end{abstract}

\section{Introduction}
The Inventory Routing Problem (IRP) integrates inventory management, vehicle routing, and delivery scheduling decisions, whereby delivery quantities are at the discretion of the distributor and are to be co-optimized with delivery routes.
The IRP most often arises in the context of \textit{vendor-managed inventory}~\citep{yao2007supply}, where a supplier anticipates inventory levels at customers and makes the replenishment decisions for products delivered to customers.
This is often a win-win strategy for both suppliers and customers: suppliers can coordinate shipment made to customers so as to save distribution costs, while customers can benefit from avoiding efforts for explicit inventory control.
The integration of inventory control and delivery planning is ubiquitously practiced in industry.
For example, the industrial gases company Air Liquide has recently pointed out that its business model is transitioning from customer-managed inventory to vendor-managed inventory, which necessitates inventory routing optimization techniques~\citep{roadef2018airliquide}. 
To this end, the ROADEF/EURO Challenge 2016~\citep{airroadef} was based on the real-life experience of Air Liquide.
Indeed, developing optimization techniques for inventory routing is of great interest for both academic research and industrial applications. 

There are many IRP variants and they can be distinguished based on different categories. 
The first category is the inventory replenishment policy. 
Replenishment policies define pre-established rules to be imposed on the quantity delivered to a customer in each visit, with the \textit{order-up-to-level} and the \textit{maximum-level} policies being most common in routing applications. The former calls for filling the inventory to the tank capacity whenever a customer is replenished~\citep{archetti2007branch}, while the latter allows for flexible replenishment and only requires that tank capacities are respected~\citep{desaulniers2015branch,avella2017single}.
Clearly, the maximum-level policy provides more flexibility for replenishment decisions and is thus more preferable from a cost-saving perspective.
The second category is the objective function.
Traditionally, the objective of an IRP is to minimize the total cost, including inventory cost and routing cost~\citep{archetti2007branch, desaulniers2015branch, avella2017single}.
Recently, some research efforts~\citep{archetti2017minimizing, baldacci2020optimal, singh2015incremental, kheiri2020heuristic, absi2020heuristic, he2020matheuristic, su2020matheuristic} have been focused on minimizing the so-called \enquote{logistic ratio}, which is the ratio of the total cost (including both inventory costs and transportation costs) to the amount of product distributed. 
Assessing the efficiency of a distribution policy through measures that consider both the distribution cost and the amount of product distributed is in fact a common industrial practice, but unfortunately, choosing the ratio function as an objective results in a fractional programming model~\citep{radzik2013fractional}, which is usually much harder to solve. 
Interested readers are referred to the work of~\citet{coelho2013thirty} for a comprehensive review on the IRP and its variants. 

A new dimension to categorize the IRP is whether to discretize the planning horizon. 
The discrete-time IRP first discretizes a given planning horizon into multiple time periods. It then assumes that customers receive their deliveries at the beginning of each period and can use them to fulfill their demands in that period~\citep{desaulniers2015branch}, irrespectively of when the actual deliveries take place as a consequence of the specific routing decisions made in each period. In this case, one only enforces inventory constraints at the beginning/end of each time period.
In contrast, the continuous-time IRP (CTIRP) explicitly accounts for the time it takes to distribute product by vehicles and manages inventory in continuous time, requiring that inventory constraints are satisfied at any time point throughout the planning horizon~\citep{lagos2020continuous}.
The practical significance of this depends of course on the exact application, but it can be argued that it is imperative in settings where feasibility of material storage and timeliness of material transfer are critical. 
This can be the case in the storage of fuels, chemicals, hazardous materials and regulated substances, among other examples. 
In these cases, vehicle arrival/departure times at customers cannot be neglected but should be explicitly accounted for by the model, in order for the latter to retain awareness of inventory levels at all times and lead to decisions that enforce them continuously.

The IRP is notoriously hard because three types of decisions have to be made simultaneously: (i)~when to serve a given customer, (ii)~how much to deliver to a given custoemr in light of other deliveries, and (iii)~how to combine customers into vehicle routes.
In the literature, the discrete-time IRP has been extensively studied. 
The work of~\citet{coelho2012consistency} considered the multi-vehicle IRP with and without consistency requirements, and applied an adaptive large neighborhood search scheme in which some subproblems are solved exactly as mixed-integer linear programs.
The first exact algorithm for IRP is from the work of~\citet{archetti2007branch} in which the authors considered a single-vehicle IRP and proposed a branch-and-cut approach that could solve to optimality benchmark instances with up to $50$ customers and $3$ time periods.
The work of~\citet{coelho2013exact} and~\citet{coelho2014improved} extended the branch-and-cut approach from~\citet{archetti2007branch} to the multi-vehicle IRP and enhanced it with new valid inequalities.
Their algorithm obtained improved upper and lower bounds for many IRP instances.
One of the most successful exact approaches for solving the discrete-time IRP was presented in~\citet{desaulniers2015branch}, where the authors proposed a sophisticated branch-price-and-cut algorithm, solving benchmark instances with up to $50$ customers and up to $5$ vehicles optimally.
In more recent work~\citet{avella2017single}, single-period tightening cuts called the \enquote{disjoint route inequalities} were proposed and incorporated into a branch-and-cut algorithm, resulting in a significant reduction in optimality gaps. 

Despite its significance, the CTIRP has received very little attention in the literature. 
The work of~\citet{dong2017solution} considered the IRP with driver rest constraints and decomposed it into an upper-level routing problem and a lower-level continuous-time scheduling problem. 
An iterative approach based on the upper and lower levels was then proposed to identify the optimal routing and scheduling design. 
\citet{fokkema2020continuous} focused on the supply-driven cyclic CTIRP in which inventory is held in containers that act as both storage and a movable transport unit, in the application of biogas transportation. 
The authors proposed a tour-based formulation that could optimally solve instances with up to $7$ customers.
The authors of~\citet{agra2022improved} proposed a location-event model and a vehicle-event model for a single-vehicle CTIRP with pickup and delivery, comparing these two models on maritime transportation instances.
Another approach was put forth in the work of~\citet{lagos2020continuous}. 
The authors first proved that, when the CTIRP instance data includes only rational numbers, there exists a sufficiently fine time discretization such that both routing and inventory replenishment decisions made in the optimal solution will always happen at these discretized time points.
Based on this theory, the authors then proposed a \textit{dynamic discretization discovery} algorithm: 
(i) the planning horizon is first discretized with different resolutions; 
(ii) for each discretization, after rounding up and down the adjusted travel times, two mixed-integer linear programs are formulated to produce an upper bound and a lower bound, respectively; 
(iii) if the difference between the two bounds is below a pre-specified tolerance, then the optimality of the CTIRP is achieved.
The novelty in the dynamic discretization discovery approach is that it idenifies exactly which times are needed to obtain an optimal, continuous-time solution by solving a sequence of relatively easy mixed-integer linear programs. 
The authors generated $90$ benchmark instances with up to $15$ customers and conducted extensive computational studies that solved $26$ of them to optimality.
\citet{lagos2022dynamic} applied this discretization approach for addressing a CTIRP with out-and-back routes where only a single customer is served in every vehicle route.
Clearly, the CTIRP is significantly more challenging than the discrete-time IRP, because vehicle routing and inventory management become much more intertwined.

In this work, we consider the deterministic CTIRP as defined in~\citet{lagos2020continuous}: 
(i)~the maximum-level inventory replenishment policy is adopted; 
(ii)~the objective is to minimize the total routing cost; and
(iii)~inventory levels are monitored and controlled continuously. 
This CTIRP has several complicated characteristics.
First, the planning horizon becomes a single but relatively long time period and, thus, vehicles have to be allowed to perform multiple trips.
Second, a customer's total product consumption across the whole planning horizon may exceed the tank capacity, necessitating multiple visits to replenish this customer during this horizon.
Third, the amount of replenishment at each customer is at the distributor's discretion, as long as the tank capacity constraints are respected, and thus allowed to vary across different visits.
Last, both vehicle routing designs and delivery scheduling decisions are now constraining the inventory control space, since they are all taken into account explicitly in continuous time.

Although the branch-price-and-cut approach proposed by~\citet{desaulniers2015branch} was quite successful in solving the discrete-time IRP, it cannot be easily extended to the CTIRP due to complications from having to model features like continuous arrival/departure times, multiple vehicle trips and multiple customer visits, among others.
Even formulating the CTIRP as a mixed-integer linear program is not a straightforward extension from the discrete-time case. 
To bridge this gap, we propose in this work a novel mixed-integer linear programming (MILP) formulation to model the continuous-time inventory routing problem. This formulation is based on the concept of duplicating customer nodes and arcs to represent routing features like multiple visits at a customer (where each node copy corresponds to a separate visit) and multiple trips by a vehicle (where each arc copy corresponds to traversal during a seaparate trip), respectively. We model the vehicle arrival/departure times at customers as continuous variables and manage stock levels in continuous time via appropriately defined inventory balance equations.
We also adapt the well-known \textit{rounded capacity inequalities} to strengthen linear programming relaxations in the context of a branch-and-cut solution process.
To elucidate the computational competitiveness of our proposed approach compared to the previous state-of-the-art, we conduct extensive studies based on $90$ benchmark instances from the literature. 
Our algorithm could solve $54$ of them to optimality within a reasonable amount of time and return a small residual gap for the remaining ones. 
We further evaluate our algorithm on newly generated benchmark instances that are inspired by realistic data from the ROADEF/EURO Challenge 2016~\citep{airroadef}. 
Out of $63$ such instances, our branch-and-cut algorithm solved $47$ of them to optimality, including a few $20$-customer instances that constitute the largest-sized CTIRP instances that have been solved to proven optimality in the open literature to-date.

The remainder of the paper is organized as follows. 
A formal problem definition is given in Section~\ref{sec:problem_definition}. 
In Section~\ref{sec:mathematical_modeling}, we present the novel MILP formulation to model the CTIRP, and then in Section~\ref{sec:tightening_techniques} we present opportunities to tighten this model and expedite its solution process. 
Section~\ref{sec:computational_studies} presents our computational studies in detail, while we conclude our work with some final remarks in Section~\ref{sec:conclusions}.

\section{Problem Definition} \label{sec:problem_definition}
The CTIRP is defined over a planning horizon $H \in \mathbb{R}_{>0}$ and on a directed graph $G = (V, E)$, where $V := \{0\} \cup V_c$ denotes the set of nodes to include a set of customers $V_c := \{1, 2, ..., n\}$ and the depot (node~$0$), and $E := \{(i, j): i \in V, j \in V \setminus \{i\} \}$ is the set of arcs. 
We assume that the depot has unlimited supply of product.
Each customer $i \in V_c$ consumes said product at a constant rate $r_i$. 
Additionally, associated with each customer $i \in V_c$ are product inventory levels $I_i^0, I_i^{\ell}, I_i^u, I_i^H$ to denote respectively the inventory at the beginning of the planning horizon, the lower and upper bounds that should be satisfied throughout the planning horizon, and the lowest acceptable inventory level at the end of the planning horizon. Without loss of generality, we assume that each customer has to be visited at least once during the planning horizon, i.e., $I_i^H + r_iH - I_i^0 > 0$ for all $i \in V_c$. 
A fleet of $K$ identical vehicles of capacity $Q$, initially located at the depot, will be used to deliver product to customers so that inventory levels at customers are maintained within the desired intervals.
Vehicles are allowed to perform multiple trips and customers are allowed to be visited multiple times during the planning horizon.
Loading at the depot and unloading at customers are considered to be instantaneous. 
However, a vehicle is allowed to remain at a customer location to dispense additional product before proceeding with the next leg of its trip.
	
Associated with each arc $(i, j) \in E$ are the time $t_{ij}$ and the cost $c_{ij}$ for a vehicle to traverse the arc. Here, we assume that both the travel time and travel cost matrices satisfy the triangle inequality.
Furthermore, we do not allow multiple vehicles to visit a customer at the same time. 
Hence, while a vehicle remains at a customer to perform additional fulfillment, no other vehicle is allowed to visit this customer (other vehicles must wait until the customer becomes available).

We do not consider inventory holding costs at the depot or at customers; hence, the total cost of a route depends only on the vehicle travel cost.
The problem's objective is then to identify the minimum-cost routing plan and schedule such that
(i)~inventory constraints are satisfied; 
(ii)~vehicle capacities are respected in each trip;
(iii)~every used vehicle returns to the depot by the end of the planning horizon.  

\section{Mathematical Modeling} \label{sec:mathematical_modeling}
In this section, we present a mathematical formulation for the CTIRP.
We highlight that the CTIRP setting deviates from archetypal vehicle routing problems in a number of non-straightforward ways, posing certain challenges for the modeler. 
For example, in the CTIRP:
(i)~we need to handle multiple visits at a customer in continuous time;
(ii)~a vehicle might deliver product multiple times while on a given customer visit;
(iii)~a vehicle might have to perform multiple trips during the planning horizon; 
(iv)~we need to monitor and constrain inventory levels continuously in time.
In what follows, we discuss how to tackle these rich features in our CTIRP model.

First, we introduce some notation (see Appendix~\ref{sec:appendix} for a summary). 
During the planning horizon, every customer consumes product at a fixed rate and the total consumption by a customer may exceed the vehicle capacity.
In such a case, several deliveries must be made to the customer by one or more vehicles. 
To that end, we will assume that customer $i \in V_c$ can be visited at most $n_i$ times during the planning horizon.
Although primarily done for modeling convenience, we note that this assumption also has a practical motivation.
Customers are generally not willing to be visited too frequently within a given interval of time. 
At the same time, a distributor usually sets a minimal amount of product per delivery, and hence, this amount can be utilized to compute a valid value for $n_i$.
Let us define $N_i := \{1, 2, ..., n_i\}$ to be the set of possible visits to customer $i \in V_c$; for convenience, we also define $N_0 := \{1\}$.
Each node of the CTIRP network is thus represented by a pair $(i, \alpha)$, where $i \in V$ indicates either the depot or a customer and $\alpha \in N_i$ indicates the visit number to $i$. 
In the case of the depot, the defined \enquote{number~1} visit is taken to apply on all occassions when a vehicle returns to the depot, whether it is for replenishment or to conclude its routing assignment.
Let $((i, \alpha), (j, \beta))$ represent the arc from node $(i, \alpha)$ to node $(j, \beta)$ in the network. 
The idea of duplicating nodes to represent multiple visits at a customer has also been employed in other vehicle routing settings, such as in maritime inventory routing~\citep{agra2018robust}, periodic vehicle routing~\citep{rothenbacher2019branch}, and green vehicle routing~\citep{koyuncu2019duplicating}.

In order to handle multiple trips by a vehicle, we propose a concept called \enquote{mode} to indicate arc copies. 
Let~$D := \{0, 1\}$ denote the set of two modes. 
In mode~\enquote{0}, a vehicle is considered to traverse an arc from $(i, \alpha)$ to $(j, \beta)$ directly, while in mode~\enquote{$1$}, it starts from node $(i, \alpha)$, makes a detour to the depot $(0, 1)$ for replenishment, and then goes to node $(j, \beta)$ for delivery.
We remark that a similar idea was also applied in the work of~\citet{karaouglan2015branch}, in which the authors duplicated the arcs so as to represent the detour trip to the depot while solving the \textit{multi-trip vehicle routing problem}.
Let $A^d$ represent the set of valid arcs in mode $d \in D$.
In particular, $A^0  :=  \{((i, \alpha), (j, \beta)): i \in V, \alpha \in N_i, j \in V \setminus \{i\}, \beta \in N_j \}$ and 
$A^1:=  \{ ((i, \alpha), (j, \beta)): i,j \in V_c, \alpha \in N_i, \beta \in N_j  \}  \setminus \{ ((i, \alpha), (i, \beta)): i \in V_c, \alpha, \beta \in N_i, \alpha \geq \beta  \}$.
Note how, in mode~\enquote{1}, the arcs $(i, \alpha, i, \beta)$ with $\alpha < \beta$ are valid, representing cases where a vehicle departs from a customer, returns to the depot for replenishment, and then visits the same customer for delivery again.
Let also $\delta_d^+(j, \beta) := \{ (i, \alpha): ((i, \alpha), (j, \beta)) \in A^d \}$ and $\delta_d^-(j, \beta) := \{(i, \alpha): ((j, \beta), (i, \alpha)) \in A^d \}$ denote the set of nodes in the network that are connected with node $(j, \beta)$ by in-coming and out-going arcs in mode $d \in D$, respectively.

\noindent \textbf{Objective Function.} The objective is to minimize the total travel cost.
The cost of traversing an arc from node $(i, \alpha)$ to node $(j, \beta)$ by a vehicle in mode $d$ is denoted by $c_{ij}^d$. Clearly, $c_{ij}^0 := c_{ij}$ (direct travel), while $c_{ij}^1 := c_{i0} + c_{0j}$ because the vehicle makes a detour to the depot.   
Let $x^d_{i \alpha j \beta}$ be a binary variable that is equal to $1$, if the arc from $(i, \alpha)$ to $(j, \beta)$ is traversed in mode $d$, and $0$ otherwise.
The objective function is thus defined as
\begin{align}
 & \min \sum_{d \in D}^{} \sum_{((i, \alpha), (j, \beta)) \in A^d}^{} c^d_{ij} x_{i\alpha j \beta}^d.   \label{eq:objective}
\end{align}

\noindent \textbf{Degree Constraints.}
Let $y_{j\beta}$ be a binary variable that is equal to $1$, if node $(j, \beta)$ is visited during the planning horizon, and $0$ otherwise. 
If a visit is made to node $(j, \beta)$, there must be an arc entering this node and an arc leaving it, as shown by~(\ref{eq:degree_in}) and~(\ref{eq:degree_out}).
\begin{align}
	& x_{i \alpha j\beta }^d \in \{0, 1\}        & \forall ((i, \alpha), (j, \beta)) \in A^d,\,\forall d \in D   \label{eq:define_x}  \\
	& y_{j \beta} \in \{0, 1\}                   & \forall \beta \in N_j,\,\forall j \in V_c  \\
	&  \sum_{d \in D}^{}  \sum_{(i, \alpha) \in \delta_d^+(j, \beta) }^{} x^d_{i \alpha j \beta} = y_{j \beta}     & \forall \beta \in N_j,\,\forall j \in V_c \label{eq:degree_in} \\ 
	&  \sum_{d \in D}^{}  \sum_{(i, \alpha) \in \delta_d^-(j, \beta) }^{} x^d_{j \beta i \alpha } = y_{j \beta}     & \forall \beta \in N_j,\,\forall j \in V_c  \label{eq:degree_out} 
\end{align}

\noindent \textbf{Fleet Size.} 
The number of used vehicles is bounded from above by $K$, the size of the fleet (initially located at the depot), as per constraint~(\ref{eq:fleet_size}). Note how the number of used vehicles can be expressed as the number of arcs that leave the depot, i.e., node $(0, 1)$, as there is no detour that starts from the depot.
\begin{align}
      &   \sum_{(j, \beta) \in \delta_0^-(0, 1) }^{}  x^0_{01j\beta} \leq K &   \label{eq:fleet_size} 
\end{align}

\noindent \textbf{Route Timing.} 
Let $[W^{\ell}_{j\beta}, W^u_{j\beta}]$ denote the time window during which a vehicle may arrive at each customer $j \in V_c$ for the $\beta^\text{th}$ visit, where $\beta \in N_j$. 
For now, one may trivially consider the interval between the earliest and latest possible arrival times, independently of the visit number $\beta$, i.e., $[W^{\ell}_{j\beta}, W^u_{j\beta}] := [t_{0j}, H - t_{j0}]$ for all $\beta \in N_j$, noting that tighter values will be proposed later~(see Section~\ref{sec:tightening_techniques}). Obviously, any applicable time windows provided directly as input data can be used instead, if tighter.

Let non-negative variables $a_{j\beta}$ and $d_{j\beta}$ respectively denote the arrival and departure time at node $(j, \beta)$, if it is visited.
Let also $\tilde{a}_{i\alpha j\beta}$ be an auxiliary non-negative variable to represent the arrival time from node $(i, \alpha)$ to node $(j, \beta)$, if this arc is traversed in any mode, and $0$ otherwise.
Constraints~(\ref{eq:define_tilde_a_00jbeta})--(\ref{eq:define_tilde_a_ialphajbeta}) then enforce that any vehicle can only leave the depot after time $0$ and has to return back to the depot by the end of the planning horizon, $H$. 
Constraints~(\ref{eq:define_arrival_departure_time}) build the relationship between variables $\tilde{a}_{i \alpha j \beta}$, $a_{j\beta}$ and $d_{j\beta}$.
Note that, since the unloading process is instantaneous, the departure from node $(j, \beta)$ can occur as soon as the arrival itself.
The time for traversing an arc from node $(i, \alpha)$ to node $(j, \beta)$ by a vehicle in mode $d$ is denoted by $T_{ij}^d$.
Clearly, $T^0_{ij} := t_{ij}$ and $T^1_{ij} := t_{i0} + t_{0j}$. 
Constraints~(\ref{eq:time_flow_constraint}) relate the departure time at node $(i, \alpha)$ to the arrival time at some node $(j, \beta)$ when a vehicle travels from $(i, \alpha)$ to $(j, \beta)$ in either mode $0$ or $1$. 
Note how these constraints are in the form of inequalities, allowing for the possibility for a vehicle to slow down while traversing an arc, or equivalently, wait near the customer premises before eventually arriving and/or right after departing.
We also remark that, even though loading and unloading products are considered to be instantaneous, fixed service times for these operations can be easily taken into account by modifying constraints~(\ref{eq:define_arrival_departure_time}) and~(\ref{eq:time_flow_constraint}) accordingly.
\begin{align}
	&a_{j\beta} \in \mathbb{R}_+  &   \beta \in N_j,\forall j \in V_c  &\label{eq:define_a} \\
	&d_{j\beta} \in \mathbb{R}_+  &   \beta \in N_j,\forall j \in V_c  &\label{eq:define_d} \\	
	&\tilde{a}_{i \alpha j \beta} \in \mathbb{R}_+  &   \forall ((i, \alpha), (j, \beta)) \in \cup_{d \in D} A^d: i, j \in V_c  &\label{eq:define_atilde} \\	
	&W^{\ell}_{j\beta} x^0_{01j\beta}  \leq  \tilde{a}_{01j \beta} \leq  W^u_{j\beta}  x^0_{01j \beta}    &   \forall (j, \beta) \in \delta_0^-(0, 1)  &\label{eq:define_tilde_a_00jbeta} \\
	& \tilde{a}_{i\alpha 01} \leq H  x^0_{i\alpha 01 }          & \forall (i, \alpha) \in \delta_0^+(0, 1)  \label{eq:define_tilde_a_jbeta00} \\
	& \tilde{a}_{i\alpha j \beta} \leq   W^u_{j\beta} \sum_{d \in D: \atop ((i, \alpha), (j, \beta)) \in A^d }^{}  x^d_{i \alpha j \beta}           & \forall ((i, \alpha), (j, \beta)) \in \cup_{d \in D} A^d: i, j \in V_c  \label{eq:define_tilde_a_ialphajbeta}  \\
	&  \sum_{(i, \alpha) \in  \cup_{d \in D} \delta_d^+(j, \beta)}^{}  \tilde{a}_{i \alpha j \beta} = a_{j \beta}  \leq   d_{j \beta}        & \forall \beta \in N_j,\,\forall j \in V_c  &\label{eq:define_arrival_departure_time} \\
	& d_{i \alpha} + \sum_{d \in D}^{} \sum_{(j, \beta) \in \delta_d^-(i, \alpha) }^{}  T^d_{ij} x^d_{i \alpha j \beta}  \leq  \sum_{(j, \beta) \in \cup_{d \in D} \delta_d^-(i, \alpha) }^{}  \tilde{a}_{i \alpha j \beta} & \forall \alpha \in N_i,\,\forall i \in V_c   \label{eq:time_flow_constraint}
\end{align}

\noindent \textbf{Non-Overlapping Visits.} When $y_{i(\alpha+1)} = 1$ and the node $(i, \alpha+1)$ is visited, that is, customer~$i$ is visited for the $(\alpha+1)^\text{th}$ time, then this customer must have been visited before, i.e., $y_{i\alpha} = 1$, as achieved by~(\ref{eq:enforce_sequence}).
Furthermore, constraints~(\ref{eq:enforce_nonoverlap_visit}) require that the recent arrival must happen after the previous departure, since a customer is not allowed to be visited by more than one vehicle at the same time.
\begin{align}
	& y_{i (\alpha + 1)} \leq y_{i \alpha}    & \forall \alpha \in N_i \setminus \{n_i\},\,\forall i \in V_c  \label{eq:enforce_sequence} \\
	& d_{i \alpha}  \leq a_{i (\alpha +1)} +  (H - t_{i0}) \left(y_{i\alpha} - y_{i \alpha+1} \right)  &  \forall \alpha \in N_i \setminus \{n_i\},\,\forall i \in V_c   \label{eq:enforce_nonoverlap_visit}
\end{align}

\noindent \textbf{Loading and Unloading Quantities.}
To model the loading and unloading quantities, we define the following continuous variables:
$f_{i \alpha j \beta}$ denotes the amount of product that a vehicle transports from node $(i, \alpha)$ to node $(j, \beta)$ in mode $0$;
$\ell_{j\beta}$ represents the amount of product that a vehicle, after finishing its current trip, reloads at the depot before immediately proceeding to node $(j, \beta)$; 
and $q_{j\beta}$ is the total amount of product delivered to customer~$j$ at the $\beta^\text{th}$ visit.
We remark that $q_{j\beta}$ might be larger than $I_i^u - I_i^{\ell}$, since vehicles are allowed to remain at a customer location and dispense product multiple times before their departure.
Constraints~(\ref{eq:flow_back_2_depot_zero}) require that any vehicle that travels from node $(i, \alpha)$ back to the depot is empty, while constraints~(\ref{eq:enforce_vehicle_capacity}) and~(\ref{eq:enforce_loading_limit}) require that the vehicle capacity is respected.
Constraints~(\ref{eq:impose unloading_upper_bound}) impose upper limits on the unloading quantities, while constraints~(\ref{eq:flow_conservation}) are the flow conservation constraints at node $(j, \beta)$, ensuring that the amount of product that flows into node $(j, \beta)$ equals the delivered amount plus the amount that flows out.
\begin{align}
	&f_{i \alpha j \beta} \in \mathbb{R}_+  &   \forall ((i, \alpha), (j, \beta)) \in A^0: j \in V_c  &\label{eq:define_f} \\	
	&\ell_{j\beta} \in \mathbb{R}_+  &   \beta \in N_j,\forall j \in V_c  &\label{eq:define_l} \\
	&q_{j\beta} \in \mathbb{R}_+  &   \beta \in N_j,\forall j \in V_c  &\label{eq:define_q} \\	
	&  f_{i \alpha 0 1} = 0                                    & \forall (i, \alpha) \in \delta_0^+(0, 1)  \label{eq:flow_back_2_depot_zero}   \\
	& f_{i \alpha j \beta} \leq Q x^0_{i \alpha j \beta}       & \forall ((i, \alpha), (j, \beta)) \in A^0: j \in V_c  \label{eq:enforce_vehicle_capacity}  \\
	& \ell_{j\beta} \leq Q  \sum_{(i, \alpha) \in \delta^+_1(j, \beta) }^{} x^1_{i \alpha j \beta }     & \forall \beta \in N_j,\,\forall j \in V_c   \label{eq:enforce_loading_limit} \\
	& q_{j \beta} \leq Q  y_{j \beta}  & \forall \beta \in N_j,\,\forall j \in V_c   \label{eq:impose unloading_upper_bound} \\
	& \sum_{(i, \alpha) \in \delta_0^+(j, \beta) }^{} f_{i \alpha j \beta} + \ell_{j \beta} = q_{j \beta} + \sum_{(i, \alpha) \in \delta_0^-(j, \beta) }^{} f_{j \beta i \alpha}     & \forall \beta \in N_j,\,\forall j \in V_c   \label{eq:flow_conservation}
\end{align}

\noindent \textbf{Inventory Limits.} Inventory limits have to be respected for each customer at any moment in the continuous planning horizon.
In order to achieve this, we specifically need to guarantee that, at the moment $a_{i\alpha}$ when a vehicle is arriving at customer $i \in V_c$, the inventory level is above the minimum desired level, $I_i^{\ell}$, and that at the moment $d_{i\alpha}$ when a vehicle is departing from customer~$i$, the inventory level is below the maximum level $I_i^u$, as enforced by constraints~(\ref{eq:inventory_lower_bound}) and~(\ref{eq:inventory_upper_bound}). In particular, note how the sums in the left-hand sides of these constraints represent the cumulative product delivery before and after the $\alpha^\text{th}$ visit, respectively.
We also remark that, during the time interval $[a_{i\alpha}, d_{i\alpha}]$, the vehicle remains at customer~$i$ and can perform multiple deliveries at zero extra cost, and thus the inventory constraints will constantly be respected in that timeframe, as long as they are satisfied at both $a_{i \alpha}$ and $d_{i \alpha}$.
At the end of the planning horizon, the inventory level at customer $i \in V_c$ should be above the target level $I_i^H$, as shown by~(\ref{eq:minimum_inventory_level}).
\begin{align}
& I_i^0y_{i\alpha} + \sum_{\alpha^{\prime} \in N_i: \atop \alpha^{\prime} < \alpha}^{} q_{i \alpha^{\prime}}  - r_i a_{i \alpha}  \geq I^{\ell}_iy_{i\alpha}  & \forall \alpha \in N_i,\,\forall i \in V_c  \label{eq:inventory_lower_bound}  \\
& I_i^0y_{i\alpha} + \sum_{\alpha^{\prime} \in N_i: \atop \alpha^{\prime} \leq \alpha}^{} q_{i \alpha^{\prime}}  - r_i d_{i \alpha}  \leq I^u_iy_{i\alpha} + \sum_{\alpha^{\prime} \in N_i: \atop \alpha^{\prime} < \alpha}^{} Q \left( y_{i\alpha^{\prime}} - y_{i\alpha} \right)  & \forall \alpha \in N_i,\,\forall i \in V_c \label{eq:inventory_upper_bound}  \\
& I_i^0 + \sum_{\alpha \in N_i}^{} q_{i \alpha}  - r_iH \geq I_i^H  & \forall i \in V_c    \label{eq:minimum_inventory_level}
\end{align}

The CTIRP aims to minimize the objective~(\ref{eq:objective}), subject to constraints~(\ref{eq:define_x})--(\ref{eq:minimum_inventory_level}).
We highlight that this is a compact MILP formulation that can handle the complex routing and scheduling features of the CTIRP, such as multiple customer visits, multiple trips by a vehicle, and continuous-time inventory management. 

\section{Model Tightening} \label{sec:tightening_techniques}
In this section, we derive various tightening inequalities to strengthen the LP relaxations of our proposed model.
In particular, we tighten time windows of any visit to a customer and use them to eliminate variables from our proposed MILP formulation.
We then consider the minimum quantity of product delivered to a customer and propose capacity constraints as strengthening cuts.

\subsection{Tightening of Time Windows}
We can tighten the time windows $[W^{\ell}_{j\beta}, W^u_{j\beta}]$ by taking advantage of the inventory restrictions.
For example, the latest time point at which node $(j, \beta)$ can be visited is when the inventory level at customer~$j$, after a postulated $\beta-1$ previous full truckload deliveries, reaches the minimum required level $I^{\ell}_j$, namely the time $(I_j^0  + (\beta-1)Q  - I_j^{\ell}) / r_j$.
Using similar reasoning, the node $(j, \beta)$ cannot be visited too early, such that even if the inventory tops off during this visit and even if an additional $n_j  - \beta$ full truckload deliveries are performed at a later time, the final inventory level is still above $I_i^H$ at the end of the planning horizon.
Hence, we may use the following tight values:
\begin{alignat} {1}
	& W^{\ell}_{j \beta} := \max \left\{ t_{0j}, H - \frac{I_j^u + \left(n_j - \beta\right)Q - I_j^H }{r_j} \right\},  \nonumber \\ 
	& W^u_{j\beta}       := \min \left\{ H - t_{j0},       \frac{ I_j^0  + \left(\beta-1\right)Q  - I_j^{\ell}}{r_j} \right\}.  \nonumber
\end{alignat}

We remark that one may further tighten these values by applying classical data preprocessing procedures as proposed in~\citet{kallehauge2007path}, noting though that we did not find this to be very beneficial in our computational studies.

\subsection{Variable Elimination}
Given node $(i, \alpha)$, node $(j, \beta)$, and mode $d \in D$, whenever $W^{\ell}_{i\alpha} + T^d_{ij} > W^u_{j\beta}$, then we can fix $x^d_{i\alpha j\beta} = 0$.
This represents a case where a vehicle cannot arrive at node $(j, \beta)$ by the end of the applicable time window even if it departs node $(i, \alpha)$ at the earliest possible time.

Moreover, let us define $\vartheta_i := I_j^H + r_jH - I_j^0$ for each customer $i \in V_c$ to signify the minimum amount of product that the distributor has to ship to customer~$i$ during the planning horizon. Note how this amount barely guarantees that the inventory level is above $I_i^{H}$ at the end of this horizon. 
We denote by $m_i$ the minimum number of visits at customer~$i$ that have to be performed. Clearly, $m_i = \ceil[\big]{\vartheta_i / Q}$ since at most a full truckload (quantity $Q$) can be delivered to customer~$i$ in each visit. 
Hence, constraints~(\ref{eq:enforce_y_to_1}) are valid.
Then, for such nodes $(i, \alpha)$ that have to be visited in every feasible solution, we can enforce valid lower bounds for arrival time variables $a_{i \alpha}$, as shown by constraints~(\ref{eq:define_arrival_time_lb}).
\begin{align}
	& y_{i \alpha } = 1                        & \forall \alpha \in N_i: \alpha \leq m_i,\,\forall i \in V_c   \label{eq:enforce_y_to_1}  \\
	&  a_{i \alpha} \geq W^{\ell}_{i \alpha}   & \forall \alpha \in N_i: \alpha \leq m_i,\,\forall i \in V_c    \label{eq:define_arrival_time_lb}
\end{align}

\subsection{Rounded Capacity Inequalities} 
The \textit{rounded capacity inequalities}~\citep{laporte1983branch} are a well-known family of strengthening inequalities that have demonstrated great effectiveness in tightening routing models and expediting their solution process within a branch-and-cut framework~\citep{toth2002vehicle}, among others.
We adapt these inequalities to the context of the CTIRP, yielding constraints
\begin{alignat}{4}
&    && \sum_{j \in S}^{} \sum_{\beta \in N_j}^{}  \left(  \sum_{(i, \alpha) \in \delta_0^+(j, \beta): \atop i \notin S}^{}  x^0_{i \alpha j\beta}  + \sum_{(i, \alpha) \in \delta_1^+(j, \beta)}^{} x^1_{i \alpha j \beta} \right) \geq k(S)  &&& \quad \forall S \subseteq V_c,  \label{eq:rounded_capacity_inequalities}
\end{alignat}
where $S \subseteq V_c$ denotes any subset of customers and $k(S)$ denotes the minimum number of vehicle visits required to serve the customers in $S$.

Note how the left-hand side in~(\ref{eq:rounded_capacity_inequalities}) represents the number of visits that are going to deliver product to customers in the set $S$, with the second term accounting for the case when vehicles return back to the depot for replenishment and then immediately serve customers in the set $S$.
In regards to the right-hand side, $k(S)$ corresponds to the optimal value of the \textit{bin-packing problem with item fragmentation}~\citep{menakerman2001bin}~(BPPIF) with bin capacity $Q$, item sizes given by the minimum demands (i.e., $\vartheta_i$) of the customers in $S$, and each item fragmented into at most $n_i$ pieces.
Recognizing that computing this optimal value exactly is as hard as the BPPIF and is therefore strongly $\mathcal{NP}$-hard, the inequality remains valid if one replaces it with the obvious lower bound $k(S) := \ceil{ \sum_{i \in S}^{}\vartheta_i / Q } $, which yields the so-called \textit{rounded capacity inequalities} (RCI).
We remark that capacity constraints similar to~(\ref{eq:rounded_capacity_inequalities}) are also used in the \textit{vehicle routing problem with split delivery}~\citep{belenguer2000lower, bianchessi2019branch}.

In the spirit of a branch-and-cut algorithm, the polynomially-sized mixed-integer linear program~(\ref{eq:objective})--(\ref{eq:define_arrival_time_lb}) and $(\ref{eq:rci_special_case})$ shall be initially presented to an MILP solver, while the RCI~(\ref{eq:rounded_capacity_inequalities}) will be dynamically separated and added into the model at each branch-and-bound node.
In order to separate the above RCI, we first construct the support graph $\bar{G} = (V, E)$ that corresponds to the LP fractional solution at a given branch-and-bound node.
In particular, we loop through all positive $\bar{x}^d_{i \alpha j \beta}$ values for $((i, \alpha), (j, \beta)) \in A^d, d \in D$, and we increment the weight of arc $(i, j) \in E$ by $\bar{x}^d_{i \alpha j \beta}$, if $d = 0$, while we increment the weights of both arcs $(i, 0) \in E$ and $(0, j) \in E$ by $\bar{x}^d_{i \alpha j \beta}$, if $d = 1$.
We emphasize that the resulting graph $\bar{G}$ is not the support graph that arises from the classic \textit{capacitated vehicle routing problem} because the degree of each customer node is not necessarily exactly $2$.
Hence, the commonly used CVRPSEP package~\citep{lysgaard2003cvrpsep} is not applicable here.
Following a separation protocol similar to the one suggested in~\citet{lysgaard2004new}, we implemented three heuristics to separate our RCI, as follows.
We first identify all \textit{connected components} in the graph $\bar{G}$ as candidates. Then, we separate the weaker \textit{fractional capacity inequalities} via the max-flow algorithm. If both heuristics fail, we then resort to a custom-built variant of the \textit{tabu search} method proposed by~\citet{augerat1998separating} to identify violated RCI.
We highlight that RCIs do not eliminate any feasible solutions since they are strengthening inequalities implied by the CTIRP model.

A special case of~(\ref{eq:rounded_capacity_inequalities}) is when $S = V_c$, as shown by constraint~(\ref{eq:rci_special_case}).
We notice that adding this constraint generally expedites the branch-and-cut algorithm; thus, in our implementation, we always append this one to our formulation from the onset. 
\begin{align}
	&    \sum_{j \in V_c}^{} \sum_{\beta \in N_j}^{} \left(x^0_{01j\beta} + \sum_{(i, \alpha) \in \delta^+_1(j, \beta)}^{} x^1_{i \alpha j \beta} \right)  \geq  \ceil[\Bigg]{ \frac{ \sum_{i \in V_c}^{}\vartheta_i }{Q} }  &  \label{eq:rci_special_case}
\end{align}

\section{Computational Studies} \label{sec:computational_studies}
In this section, we test our branch-and-cut algorithm on CTIRP benchmark instances and compare it against the state-of-the-art approach from~\citet{lagos2020continuous}.
Our algorithm was implemented in C++ and the mixed-integer linear program was solved using the Gurobi Optimizer~$11.0.0$ through the C application programming interface.
User cuts were implemented via the Gurobi callback function, while branching prioritized binary variables $y_{i\alpha}$ over variables $x^{d}_{i\alpha j \beta}$ via the Gurobi attribute-setting function.
All solver settings were kept at their default, except the optimality gap tolerances, with the relative tolerance set to $0$ (effectively disabled) and the absolute one to the value of $0.0099$. Note that the latter suffices for any optimal solutions to be exact, since all entries of the travel cost matrix in all datasets used have at most two decimal digits and, hence, so does the objective value of any feasible solution.
All computations were performed on a server with dual Intel~Xeon~E5-2689v4 CPUs clocking at 3.10~GHz. 
The 128~GB of available RAM was shared among $10$ copies of the algorithm running in parallel on this machine. 
Each instance was solved by one copy of the algorithm limited to use a single thread.

\subsection{Benchmark Instances}
We consider two CTIRP datasets for testing our algorithm.
The first one stems from the work of~\citet{lagos2020continuous}, in which the authors generated two types of CTIRP instances, namely \texttt{clustered}~(C) and \texttt{random}~(R).
There are $45$ instances for each type, with the number of customers ranging from $5$~to~$15$ and the number of vehicle ranging from $3$~to~$15$. 
All instances can be found at \url{https://github.com/felipelagos/cirplib}.
As per~\citet{lagos2020continuous}, the travel costs and travel times are taken as the Euclidean distances rounded to two decimal places.

Furthermore, to evaluate the performance of our algorithm on solving realistic CTIRP data, we generate another dataset inspired by the ROADEF/EURO Challenge~2016 competition~\citep{airroadef}, which compiled data from the practice of distributing liquefied gases to customers in a vendor-managed inventory setting.
In particular, we obtain the \enquote{version $1.1$ set A} dataset from this competition and apply the following modifications: 
(i) the fleet is reduced to a homogeneous one; 
(ii) the driver-trailer scheduling aspect is neglected; 
(iii) the product consumption rate is considered to be constant throughout the planning horizon; 
(iv) the objective is to minimize the travel cost, rather than the logistic ratio.
We highlight that this dataset includes $9$ medium-size instances stemming from real-life practice, with the number of customers ranging from $53$~to~$89$.
To create our CTIRP instances, we choose the first $5$, $7$, $10$, $12$, $15$, $17$ and~$20$ customers from each original instance, thus resulting in a suite of $63$ instances that we shall refer to as the \texttt{roadef} (RF) type. 
The number of vehicles in our instances ranges from $1$~to~$6$. 
Our newly generated instances are named \enquote{RF-X-nY-kZ} to signify how the X$^{th}$ original instance was adapted to include Y~customers and Z~vehicles.
All instance data can be found at \url{http://gounaris.cheme.cmu.edu/datasets/cirp/}.

In all our experiments, we limit the number of visits to each customer $i \in V_c$ to two plus the minimal number of visits at customer~$i$, i.e., $n_i := m_i + 2$. Whereas any selection of $n_i$ might be limiting for obtaining a strictly better solution, we highlight that, for no instance in our computational studies, did we notice an improvement in the optimal objective value when we increased this number to $m_i+3$. In this manner, we have empirically determined the above choice of $n_i$ to be satisfactory.

\subsection{Computational Results on Literature Benchmarks}  \label{sec:comparison_with_existing_method}
We first test our branch-and-cut algorithm (denoted by \enquote{Branch-and-Cut}) on the first dataset (\texttt{clustered} and \texttt{random} instances) and compare its performance with that from the work of~\citet{lagos2020continuous}~(denoted by \texttt{LBS20}). Since that study applied both lower and upper bounding procedures (each with a corresponding time limit of $2$~hours) to obtain the best possible results, in our method we impose a time limit of $4$~hours for each instance.
In order to analyze the effect of incorporating RCI as part of the solution process, we also consider a variant of our algorithm in which RCI are disabled. In this case, no user cuts are added to the model and, hence, we denote it by \enquote{Gurobi (default)}.
We also highlight that, in the \enquote{Branch-and-Cut} approach, all Gurobi-generated, general-purpose cuts were disabled because we observed that using them decreased overall computational performance.

The computational results from each approach are compared in Table~\ref{table:comparison_against_lagos}.
Columns~\enquote{\#~cust.} and~\enquote{\#~inst.} report the number of customers in an instance and the number of instances of such an input size, respectively.
For each respective algorithm, the column \enquote{\#~opt.} reports the number of instances that were solved to optimality within the time limit, while column \enquote{\#~t.l.} reports the number of instances for which the algorithm was terminated due to the time limit but with valid lower and upper bounds identified.
For the instances solved to optimality, we report the average solution time (in seconds) and the average number of explored nodes in columns~\enquote{Avg.~t (s)} and \enquote{Avg. \# nodes}, respectively. These entries are rounded to the nearest integer.
For the unsolved instances but with feasible solutions identified, we report in column~\enquote{Avg.~gap~(\%)} the average residual gap, defined as $(UB-LB)/UB$ at the time limit.

Out of the totality of $90$ CTIRP instances, \texttt{LBS20} solved $26$ of them to optimality, returning an average residual gap of 6.45\% in the remaining~$64$.
In contrast, the proposed \enquote{Branch-and-Cut} algorithm solved $54$ instances to guaranteed optimality with an an average gap of $1.19\%$ for $35$ of the remaining instances.
In particular, our algorithm solved to optimality all but two of the instances up to $7$ customers, half of the $10$-customer instances, as well as a handful of the $12$- and $15$-customer instances. For a single \texttt{random} instance with $15$ customers, our algorithm could not identify a feasible solution within the time limit of $4$ hours.
Adressing the CTIRP using the \enquote{Gurobi (default)} approach was found to be moderately effective, being able to prove optimality for $42$ instances, while returning an average residual gap of $4.20\%$ for 45 instances. \enquote{Gurobi (default)} did not identify feasible solutions to a total of three $15$-customer instances.
Whereas incorporating RCI cuts has a noticeable effect towards improving overall tractability, it is noteworthy that the no-RCI variant can still close the gap in approximately half the instances considered. Indeed, a practitioner can readily take advantage of this latter approach's performance, since it merely calls for implementing and solving our formulation as a monolithic MILP model, without incorporating the cut separation routine.

Tables~\ref{table:branch_and_cut_clustered} and~\ref{table:branch_and_cut_random} present detailed results on the \texttt{clustered}~(C) and \texttt{random}~(R) instances, respectively, for our proposed \enquote{Branch-and-Cut} algorithm.
In these tables, if an instance could be solved to optimality within $4$ hours, then \enquote{Opt [UB]} reports the corresponding optimal objective value, while \enquote{t (s) [LB]} reports the time (in seconds, rounded to the nearest integer) to solve the instance to optimality. 
Otherwise, the columns respectively report in brackets the best upper and lower bounds found by the time limit.
We also report the average number of visits across all customers in the returned best known solution to each instance (column \enquote{Avg. \# visits}) and the number of RCI that have been dynamically added (column \enquote{\# RCI}).
Notably, the average number of visits to every customer ranges from $1.2$ to $2.2$, which is far from the maximum number of allowable visits, $n_i$, that generally ranges from $3$~to~$5$. This indicates that our choice of the number of allowable visits is most likely not restrictive. 

\begin{sidewaystable}[htbp] 
	\centering
	\caption{Comparison among \texttt{LBS20} and our proposed approaches on benchmark instances} 
	\renewcommand{\arraystretch}{0.9}
	\setlength\tabcolsep{0.1mm}
	\begin{tabular}{c c ccc ccccc ccccc}
		\toprule
		\multirow{2}{*}{\# cust.} 
		& \multirow{2}{*}{\# inst.} 
		& \multicolumn{3}{c}{\texttt{LBS20}}  
		& \multicolumn{5}{c}{Gurobi (default)}
		& \multicolumn{5}{c}{Branch-and-Cut} \\
		\cmidrule(l){3-5} 
		\cmidrule(l){6-10}   
		\cmidrule(l){11-15}  
		\parbox{0.1\textwidth}{}  
		& \parbox{0.1\textwidth}{} 
		& \parbox{0.065\textwidth}{\centering \# opt.}  
		& \parbox{0.065\textwidth}{\centering \# t.l.}  
		& \parbox{0.065\textwidth}{\centering Avg. gap (\%)}  
		& \parbox{0.065\textwidth}{\centering \# opt.}  
		& \parbox{0.065\textwidth}{\centering \# t.l.}
		& \parbox{0.065\textwidth}{\centering Avg. t (s)} 		
		& \parbox{0.065\textwidth}{\centering Avg. \# nodes} 	
		& \parbox{0.065\textwidth}{\centering Avg. gap (\%)}   
		& \parbox{0.065\textwidth}{\centering \# opt.}  
		& \parbox{0.065\textwidth}{\centering \# t.l.}
		& \parbox{0.065\textwidth}{\centering Avg. t (s)}
		& \parbox{0.065\textwidth}{\centering Avg. \# nodes} 	 		
		& \parbox{0.065\textwidth}{\centering Avg. gap (\%)}   
		\\
		\midrule
	  	   &        &      &     &        &      &       &\textbf{clustered}     \\
		5  &   9    &   7  &   2 &  1.55  &  9   &  0    &  17   &1,474   & --   &  9  &  0    & 2     & 156     & --   \\
		7  &   9    &   3  &   6 &  1.92  &  8   &  1    &  3,351& 85,557 & 0.53 &  7  &  2    & 1,385 & 10,391  & 0.58   \\
		10 &   9    &   1  &   8 &  1.90  &  1   &  8    &  34   & 1,936   & 2.50 &  3  &  6    & 2,514 & 52,285  & 1.16  \\
		12 &   9    &   1  &   8 &  3.26  &  0   &  9    &  --   & --     & 4.15 &  2  &  7    & 36    & 798     & 0.87  \\
		15 &   9    &   1  &   8 &  8.55  &  0   &  8    &  --   & --     & 6.83 &  1  &  8    & 340   & 3,268   & 1.19  \\	 
	   	   &        &      &     &        &      &       &  \textbf{random}       \\
		5  &   9    &   5  &   4 &  2.93  &  9   &  0    &  2     & 217       & --   & 9   & 0     & 1    & 57     & --    \\
		7  &   9    &   5  &   4 &  3.27  &  9   &  0    &  293   & 30,650    & --   & 9   & 0     & 27   & 622    & --   \\
		10 &   9    &   3  &   6 &  4.10  &  5   &  4    &  3,567 & 113,473   & 2.36 & 7   & 2     & 91   & 3,840  & 2.05  \\
		12 &   9    &   0  &   9 &  8.93  &  1   &  8    &  1,420 & 60,909    & 3.78 & 5   & 4     & 3,244& 52,366 & 2.03   \\
		15 &   9    &   0  &   9 &  17.67 &  0   &  7    &  --    & --        & 5.27 & 2   & 6     & 3,573& 875,006& 0.94    \\	
		\midrule
		Total &  90  & 26   & 64  &      &  42  & 45  &       &  &      & 54   &35   &       &  &  \\
		Avg.  &      &      &     & 6.45 &      &     & 1,164 &  & 4.20 &      &     & 776   &  & 1.19  \\ 
		\bottomrule
	\end{tabular}
	\label{table:comparison_against_lagos}
\end{sidewaystable}

\begin{table}[htbp] 
	\centering
	\caption{Detailed results for the Branch-and-Cut algorithm on $45$ \texttt{clustered} instances} 
	\renewcommand{\arraystretch}{0.9}
	\setlength\tabcolsep{0.6mm}
	\begin{tabular}{ccccc c ccccc}
		\toprule			
		\parbox{0.1\textwidth}{\centering Inst.}  		
		& \parbox{0.08\textwidth}{\centering Opt [UB]}  		
		& \parbox{0.08\textwidth}{\centering t (s) [LB]} 
		& \parbox{0.1\textwidth}{\centering Avg. \# visits} 
		& \parbox{0.06\textwidth}{\centering \# RCI}  
		& \hspace{0.8cm}	
		& \parbox{0.1\textwidth}{\centering Inst.} 
		& \parbox{0.08\textwidth}{\centering Opt [UB]}  		
		& \parbox{0.08\textwidth}{\centering t (s) [LB] } 
		& \parbox{0.1\textwidth}{\centering Avg. \# visits}
		& \parbox{0.06\textwidth}{\centering \# RCI}  	
		\\
		\midrule	
	C5U1Q1 &37.85 &1 &1.4 &7 &  &C10U2Q3 &[56.10] &[54.81] &1.7 &74,483\\
	C5U1Q2 &30.61 &1 &1.4 &16 &  &C10U3Q1 &104.37 &2,029 &2.1 &164\\
	C5U1Q3 &28.24 &1 &1.2 &6 &  &C10U3Q2 &[83.35] &[83.16] &1.9 &60,068\\
	C5U2Q1 &50.09 &1 &1.8 &8 &  &C10U3Q3 &[65.52] &[64.87] &1.8 &48,330\\
	C5U2Q2 &38.28 &1 &1.6 &53 &  &C12U1Q1 &85.07 &13 &1.4 &64\\
	C5U2Q3 &30.61 &1 &1.4 &12 &  &C12U1Q2 &[69.28] &[68.79] &1.7 &40,720\\
	C5U3Q1 &57.74 &14 &2.2 &41 &  &C12U1Q3 &[57.07] &[56.31] &1.5 &66,311\\
	C5U3Q2 &43.36 &1 &1.8 &12 &  &C12U2Q1 &[104.29] &[104.11] &2.0 &793\\
	C5U3Q3 &36.16 &1 &1.8 &10 &  &C12U2Q2 &[80.01] &[79.28] &1.8 &23,754\\
	C7U1Q1 &52.14 &2,442 &1.6 &3,866 &  &C12U2Q3 &[69.49] &[68.80] &1.7 &36,368\\
	C7U1Q2 &[42.99] &[42.72] &1.7 &85,149 &  &C12U3Q1 &124.43 &58 &2.1 &87\\
	C7U1Q3 &31.17 &23 &1.4 &178 &  &C12U3Q2 &[101.91] &[101.05] &2.1 &27,154\\
	C7U2Q1 &59.75 &6 &1.9 &46 &  &C12U3Q3 &[79.82] &[78.91] &1.8 &14,043\\
	C7U2Q2 &[50.48] &[50.21] &2.0 &50,632 &  &C15U1Q1 &99.78 &340 &1.4 &356\\
	C7U2Q3 &42.18 &55 &1.6 &205 &  &C15U1Q2 &[77.32] &[76.59] &1.5 &16,240\\
	C7U3Q1 &78.08 &5 &2.1 &41 &  &C15U1Q3 &[64.53] &[63.89] &1.5 &46,765\\
	C7U3Q2 &59.44 &1 &1.6 &14 &  &C15U2Q1 &[122.45] &[122.06] &1.9 &4,943\\
	C7U3Q3 &50.18 &7,165 &1.7 &10,686 &  &C15U2Q2 &[95.02] &[93.98] &1.7 &3,966\\
	C10U1Q1 &70.97 &10 &1.4 &115 &  &C15U2Q3 &[78.88] &[76.59] &1.7 &11,443\\
	C10U1Q2 &[55.62] &[54.78] &1.6 &104,290 &  &C15U3Q1 &[147.06] &[147.02] &2.1 &3,802\\
	C10U1Q3 &[45.24] &[44.84] &1.6 &101,412 &  &C15U3Q2 &[118.36] &[116.43] &2.1 &19,013\\
	C10U2Q1 &86.02 &5,504 &1.8 &139 &  &C15U3Q3 &[95.14] &[93.62] &1.9 &9,143\\
	C10U2Q2 &[65.92] &[65.22] &1.8 &69,067\\
		\bottomrule
	\end{tabular}
	\label{table:branch_and_cut_clustered}
\end{table}

\begin{table}[htbp] 
	\centering
	\caption{Detailed results for the Branch-and-Cut algorithm on $45$ \texttt{random} instances} 
	\renewcommand{\arraystretch}{0.9}
	\setlength\tabcolsep{0.6mm}
	\begin{tabular}{ccccc c ccccc}
		\toprule			
		\parbox{0.1\textwidth}{\centering Inst.}  		
		& \parbox{0.08\textwidth}{\centering Opt [UB]}  		
		& \parbox{0.08\textwidth}{\centering t (s) [LB]} 
		& \parbox{0.1\textwidth}{\centering Avg. \# visits} 
		& \parbox{0.06\textwidth}{\centering \# RCI}  
		& \hspace{0.8cm}	
		& \parbox{0.1\textwidth}{\centering Inst.} 
		& \parbox{0.08\textwidth}{\centering Opt [UB]}  		
		& \parbox{0.08\textwidth}{\centering t (s) [LB] } 
		& \parbox{0.1\textwidth}{\centering Avg. \# visits}
		& \parbox{0.06\textwidth}{\centering \# RCI}  	
		\\
		\midrule
R5U1Q1 &36.42 &1 &1.4 &2 &  &R10U2Q3 &64.02 &6 &1.6 &38\\
R5U1Q2 &29.06 &1 &1.2 &5 &  &R10U3Q1 &116.04 &35 &2.2 &56\\
R5U1Q3 &28.45 &1 &1.2 &14 &  &R10U3Q2 &[93.93] &[92.09] &1.9 &37,792\\
R5U2Q1 &41.42 &1 &1.8 &6 &  &R10U3Q3 &82.92 &245 &1.8 &438\\
R5U2Q2 &36.51 &1 &1.6 &6 &  &R12U1Q1 &[101.88] &[99.90] &1.8 &35,978\\
R5U2Q3 &30.90 &1 &1.4 &5 &  &R12U1Q2 &75.45 &19 &1.4 &137\\
R5U3Q1 &45.31 &1 &2.0 &6 &  &R12U1Q3 &66.09 &12,899 &1.2 &41,487\\
R5U3Q2 &39.44 &8 &2.0 &28 &  &R12U2Q1 &[118.75] &[118.48] &2.0 &3,007\\
R5U3Q3 &33.54 &1 &1.6 &24 &  &R12U2Q2 &98.18 &882 &1.8 &1,721\\
R7U1Q1 &63.03 &1 &1.6 &7 &  &R12U2Q3 &76.27 &2,079 &1.6 &2,776\\
R7U1Q2 &43.93 &1 &1.4 &14 &  &R12U3Q1 &139.88 &342 &2.1 &413\\
R7U1Q3 &40.43 &167 &1.6 &1,044 &  &R12U3Q2 &[111.63] &[107.60] &1.8 &22,795\\
R7U2Q1 &69.37 &1 &1.7 &7 &  &R12U3Q3 &[95.74] &[93.50] &1.6 &37,545\\
R7U2Q2 &60.12 &1 &1.7 &33 &  &R15U1Q1 &[118.48] &[117.52] &1.9 &79,728\\
R7U2Q3 &43.93 &9 &1.3 &75 &  &R15U1Q2 &87.67 &3,971 &1.6 &1,182\\
R7U3Q1 &81.58 &2 &2.0 &21 &  &R15U1Q3 &75.25 &3,175 &1.3 &10,042\\
R7U3Q2 &67.06 &4 &1.9 &27 &  &R15U2Q1 &[139.63] &[138.72] &1.9 &28,466\\
R7U3Q3 &57.47 &55 &1.9 &321 &  &R15U2Q2 &[115.99] &[114.43] &1.8 &72,226\\
R10U1Q1 &[87.72] &[85.84] &1.9 &28,759 &  &R15U2Q3 &[89.45] &[88.69] &1.6 &22,512\\
R10U1Q2 &64.02 &9 &1.6 &68 &  &R15U3Q1 &[162.61] &[161.98] &2.1 &2,858\\
R10U1Q3 &54.35 &212 &1.3 &954 &  &R15U3Q2 &[129.86] &[127.75] &1.8 &24,344\\
R10U2Q1 &99.48 &6 &1.7 &57 &  &R15U3Q3 & No sol.  &[111.07] & -- &58,630\\
R10U2Q2 &82.92 &127 &1.8 &126\\
		\bottomrule
	\end{tabular}
	\label{table:branch_and_cut_random}
\end{table}

To take a closer look at model tightness, we consider the four instances for which detailed computational results were presented in the appendix of~\cite{lagos2020continuous} and report in Table~\ref{table:comparison_against_the_literature2} relevant metrics, such as lower and upper bounds, gaps, and the overall solution times required, for each of the two approaches. Note that, for the~\texttt{LBS20} approach, we focus only on the results under the discretization parameter~($H/\Delta$) value (out of 21 different choices) that produced the smallest residual gaps in the corresponding instance.
We highlight that the \textit{root node} lower bound produced by our formulation (usually within a few seconds) is generally better than the \textit{final} lower bound from \texttt{LBS20}, indicating that our proposed model corresponds to a tighter relaxation.

\begin{table}[htbp] 
	\centering
	\caption{Detailed results for \texttt{LBS20} and our Branch-and-Cut algorithm on $4$ benchmark instances} 
	\renewcommand{\arraystretch}{0.9}
	\setlength\tabcolsep{0.1mm}
	\begin{tabular}{cc  cccc ccccc}
		\toprule
		  \multirow{2}{*}{Inst.} 
		& \multirow{2}{*}{\# cust.}   
		& \multicolumn{4}{c}{\texttt{LBS20}}  	
		& \multicolumn{5}{c}{Branch-and-Cut} \\
		  \cmidrule(l){3-6} 
		  \cmidrule(l){7-11} 
	 	  \parbox{0.1\textwidth}{}  
		& \parbox{0.1\textwidth}{}  
  		& \parbox{0.08\textwidth}{\centering $H/\Delta$}  	
		& \parbox{0.08\textwidth}{\centering UB}  		
		& \parbox{0.08\textwidth}{\centering LB}  
		& \parbox{0.08\textwidth}{\centering Gap (\%)}  		
		& \parbox{0.08\textwidth}{\centering UB}  
		& \parbox{0.08\textwidth}{\centering Root LB}  	
		& \parbox{0.08\textwidth}{\centering LB}  
		& \parbox{0.08\textwidth}{\centering Gap (\%)}  
		& \parbox{0.08\textwidth}{\centering t (s)} 			
		\\
		\midrule
		R5U2Q2   &  5 & 9 & 39.23 &  36.42   & 7.16  & 36.51 & 36.42  & 36.51 & 0.00  & 1 \\
		R10U1Q2  & 10 & 18 & 64.11 &  59.18   & 7.69  & 64.02 & 58.56  & 64.02 & 0.00  & 9 \\
		R10U2Q2  & 10 & 18 & 87.92 &  73.62   & 16.26  & 82.92 & 75.54  & 82.92 & 0.00  & 127 \\
		C15U2Q2  & 15 & 30 & 97.00 &  86.27   & 11.06  & 95.02 & 92.64   & 93.98 & 1.09  & 14,400 \\
		\bottomrule
	\end{tabular}
	\label{table:comparison_against_the_literature2}
\end{table}

\begin{figure}[!htb]
	\captionsetup[subfigure]{belowskip=0pt}
	\centering
	\begin{subfigure}[b]{0.5\textwidth}\centering
		\includegraphics[height=6.2cm]{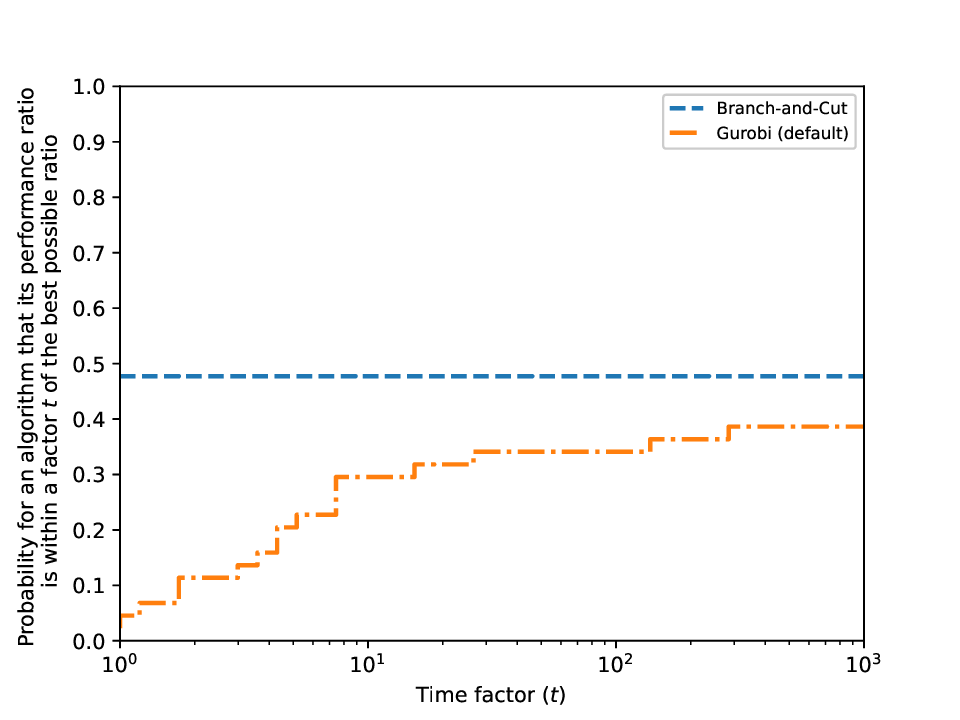}
		\caption{\texttt{clustered} instances}
		\label{figure:performance_profile_cluster}
	\end{subfigure}~%
	\begin{subfigure}[b]{0.5\textwidth}\centering
		\includegraphics[height=6.2cm]{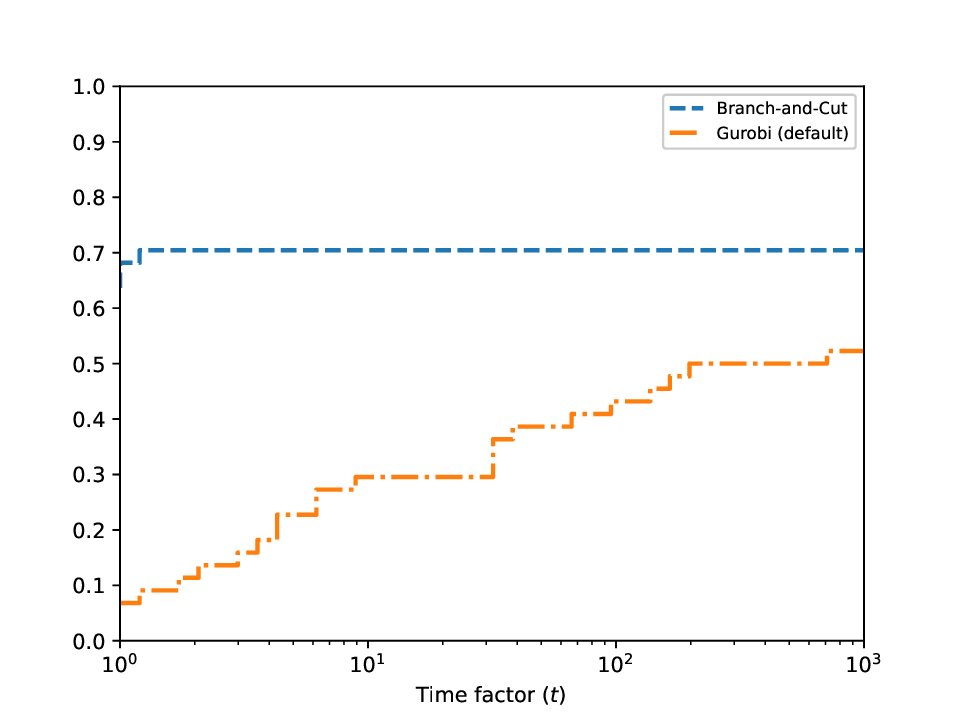}
		\caption{\texttt{random} instances}
		\label{figure:performance_profile_random}
	\end{subfigure}
	\caption{Log-scaled performance profiles across all benchmark instances. For each curve, the value at $t = 0$ indicates the fraction of benchmark instances that can be solved the fastest via the corresponding approach, while the limiting value at $t \to \infty$ indicates the fraction of instances that could be solved within the time limit.}
	\label{figure:performance_profile}
\end{figure}

We now revisit Table~\ref{table:comparison_against_lagos} for evaluating the effect of RCI on the branch-and-bound process.
Compared with the \enquote{Gurobi (default)} version, the \enquote{Branch-and-Cut} algorithm with RCI enabled as strengthening inequalities solved to optimality $12$ additional benchmark instances, including a few of the larger $15$-customer ones, while producing much smaller branch-and-bound trees.
We also note that solution times and residual gaps for unsolved instances were decreased across the board, indicating the effectiveness of RCI in expediting the branch-and-bound based search process. 
To better illustrate this, performance profiles~\citep{dolan2002benchmarking} are presented in Figure~\ref{figure:performance_profile}, which generally reveal that our \enquote{Branch-and-Cut} algorithm can solve more instances to optimality with less computational time.
Finally, it is worthy to note that the total number of RCI that were identified and introduced as tightening inequalities is in the order of hundreds for small-sized CTIRP instances, increasing to tens of thousands for the large-sized ones (see Tables~\ref{table:branch_and_cut_clustered} and~\ref{table:branch_and_cut_random}).

\subsection{Computational Results on \texttt{roadef} Instances} \label{sec:computational_results_roadef}
We now turn our attention to the set of the newly constructed \texttt{roadef} instances and evaluate the performance of the \enquote{Branch-and-Cut} algorithm on those.
The imposed time limit for this study was chosen as $10$ hours.
A synopsis of our computational results is presented in Table~\ref{table:computational_results_roadef} while the detailed results are reported in Table~\ref{table:branch_and_cut_roadef}.
Out of the $63$ instances in this dataset, our \enquote{Branch-and-Cut} algorithm solved $59$ to optimality and returned an average residual gap of $0.07\%$ for the remaining $4$ instances. 
In particular, our algorithm solved all instances with up to $12$ customers as well as the vast majority of the larger ones with up to $20$ customers. 
Notably, our approach was able to identify feasible solutions in all instances.
The solution times generally increase with problem size, ranging from minutes to a few hours.
Interestingly, the average customer is visited several times (about $3$) during the planning horizon.
Here, we should point out that the \enquote{Gurobi (default)} approach was again outperformed by our \enquote{Branch-and-Cut} approach. For reference, the former could solve to optimality only $20$ of these \texttt{roadef} instances within the same time limit. For the sake of brevity, we omit presenting the detailed computational results for that method.

\begin{table}[htbp] 
	\centering
	\caption{Results synopsis for the Branch-and-Cut algorithm on \texttt{roadef} instances} 
	\renewcommand{\arraystretch}{0.9}
	\setlength\tabcolsep{0.1mm}
	\begin{tabular}{cc cc cc}
		  \toprule
		  \parbox{0.12\textwidth}{\centering \# cust.}  
		& \parbox{0.12\textwidth}{\centering \# inst.} 
		& \parbox{0.12\textwidth}{\centering \# opt.}  
		& \parbox{0.12\textwidth}{\centering \# t.l.}  
		& \parbox{0.12\textwidth}{\centering Avg. t (s)} 		
		& \parbox{0.12\textwidth}{\centering Avg. gap (\%)}  
		\\
		\midrule
		5   &   9   & 9    &  0  &  39   &  --   \\
		7   &   9   & 9    &  0  & 388   &  --   \\
		10  &   9   & 9    &  0  & 1,605 &  --   \\
		12  &   9   & 9    &  0  & 6,480 &  --   \\
		15  &   9   & 8    &  1  & 3,143 &  0.12 \\
		17  &   9   & 8    &  1  & 5,461 &  0.04 \\
		20  &   9   & 7    &  2  & 10,434&  0.06 \\
		\midrule
	  Total &  63   & 59   & 4   &       &       \\
	   Avg. &       &      &     & 3,703 & 0.07  \\
		\bottomrule
	\end{tabular}
	\label{table:computational_results_roadef}
\end{table}

\begin{sidewaystable}[htbp] 
	\centering
	\caption{Detailed results for the Branch-and-Cut algorithm on $63$ \texttt{roadef} instances} 
	\renewcommand{\arraystretch}{0.8}
	\setlength\tabcolsep{0.6mm}
	\begin{tabular}{ccc c ccc  c ccc}
		\toprule			
		\parbox{0.1\textwidth}{\centering Inst.}  		
		& \parbox{0.08\textwidth}{\centering Opt [UB]}  		
		& \parbox{0.08\textwidth}{\centering t (s) [LB]} 
		& \hspace{0.8cm}	
		& \parbox{0.1\textwidth}{\centering Inst.} 
		& \parbox{0.08\textwidth}{\centering Opt [UB]}  		
		& \parbox{0.08\textwidth}{\centering t (s) [LB] } 
		& \hspace{0.8cm}	
		& \parbox{0.1\textwidth}{\centering Inst.} 
		& \parbox{0.08\textwidth}{\centering Opt [UB]}  		
		& \parbox{0.08\textwidth}{\centering t (s) [LB] } 	
		\\
		\midrule
RF-3-n5-k1 &1,039.5 &4 &  & RF-6-n10-k3 &1,998.9 &539 &  & RF-9-n15-k3 &2,000.7 &1,970\\
RF-4-n5-k1 &1,039.5 &4 &  & RF-7-n10-k1 &1,652.4 &4,808 &  & RF-10-n15-k2 &2,520.0 &156\\
RF-5-n5-k1 &542.4 &1 &  & RF-8-n10-k2 &1,837.5 &5,981 &  & RF-11-n15-k3 &2,520.0 &1,698\\
RF-6-n5-k1 &542.4 &1 &  & RF-9-n10-k2 &1,631.1 &55 &  & RF-3-n17-k3 &4,886.0 &5,248\\
RF-7-n5-k1 &492.0 &39 &  & RF-10-n10-k2 &1,869.6 &26 &  & RF-4-n17-k2 &4,603.2 &782\\
RF-8-n5-k1 &524.7 &143 &  & RF-11-n10-k2 &1,869.6 &950 &  & RF-5-n17-k4 &[2,305.8] &[2,304.9]\\
RF-9-n5-k1 &480.6 &98 &  & RF-3-n12-k2 &3,804.5 &2,138 &  & RF-6-n17-k5 &2,421.9 &1,864\\
RF-10-n5-k1 &1,323.9 &54 &  & RF-4-n12-k2 &3,757.6 &1,084 &  & RF-7-n17-k3 &2,094.0 &7,989\\
RF-11-n5-k1 &1,323.9 &6 &  & RF-5-n12-k4 &2,034.9 &10,853 &  & RF-8-n17-k4 &2,366.1 &5,653\\
RF-1-n7-k2 &6,309.0 &9 &  & RF-6-n12-k3 &2,157.6 &24,267 &  & RF-9-n17-k3 &2,050.8 &8,442\\
RF-2-n7-k2 &6,309.0 &5 &  & RF-7-n12-k2 &1,760.7 &10,089 &  & RF-10-n17-k4 &3,145.8 &7,582\\
RF-3-n7-k2 &2,618.7 &13 &  & RF-8-n12-k2 &1,945.8 &7,056 &  & RF-11-n17-k4 &3,173.7 &6,125\\
RF-4-n7-k1 &2,618.7 &2,122 &  & RF-9-n12-k2 &1,739.4 &748 &  & RF-3-n20-k4 &5,518.8 &22,322\\
RF-5-n7-k2 &927.3 &9 &  & RF-10-n12-k2 &1,996.5 &537 &  & RF-4-n20-k3 &7,410.2 &12,360\\
RF-6-n7-k1 &1,007.1 &11 &  & RF-11-n12-k2 &1,996.5 &1,546 &  & RF-5-n20-k4 &[2,315.4] &[2,313.9]\\
RF-7-n7-k1 &739.8 &17 &  & RF-3-n15-k2 &4,552.8 &2,979 &  & RF-6-n20-k4 &2,436.9 &10,937\\
RF-8-n7-k1 &781.5 & 205 &  & RF-4-n15-k3 &4,027.8 &830 &  & RF-7-n20-k4 &2,306.4 &5,232\\
RF-9-n7-k1 &731.7 & 1,104 &  & RF-5-n15-k4 &[2,277.3] &[2,274.6] &  & RF-8-n20-k6 &[2,583.0] &[2,581.8]\\
RF-3-n10-k2 &3,086.3 &42 &  & RF-6-n15-k4 &2,392.8 &7,979 &  & RF-9-n20-k4 &2,280.0 &2,966\\
RF-4-n10-k1 &2,913.4 &1,707 &  & RF-7-n15-k2 &2,052.6 &9,081 &  & RF-10-n20-k4 &3,375.6 &4,516\\
RF-5-n10-k3 &1,879.8 &340 &  & RF-8-n15-k3 &2,314.2 &448 &  & RF-11-n20-k4 &3,427.5 &14,708\\
		\bottomrule
	\end{tabular}
	\label{table:branch_and_cut_roadef}
\end{sidewaystable}

\section{Conclusions} \label{sec:conclusions}

In this work, we consider the continuous-time inventory routing problem (IRP), a challenging setting with ubiquitous practical application. 
This problem often arises in the context of vendor-managed inventory systems where the distributor anticipates the customer's product usage and is charged with actively managing inventory levels at the customer site. Notably, the level of the inventory must be maintained within acceptable limits in continuous time, unlike earlier IRP models that only enforced this in discrete time points during the planning horizon.

We propose a novel mixed-integer linear programming formulation that incorporates multi-trip and multi-visit features as well as continuous-time inventory management.
To expedite the solution process, we propose various types of tightening techniques, including among others the adaptation of well-known rounded capacities inequalities.
We conduct extensive computational studies on $90$ benchmark instances from the literature, revealing the effectiveness of our branch-and-cut algorithm and its outperforming of the previous state-of-the-art approach. 
In particular, we close $28$ previously open CTIRP benchmark instances and return an average gap below $1.2\%$ for the unsolved instances. 
We also further computational studies using newly constructed benchmarks based on real-life data. 
These show that our proposed algorithm could solve CTIRP instances of up to $20$ customers and $4$ vehicles.

Overall, the CTIRP constitutes an extremely challenging problem, especially as one compares with its discrete-time IRP counterpart for which state-of-the art approaches can address instances with up to $100$ customers.
It remains a question for the practitioner whether the more rigorous feasibility guarantees afforded by tracking and enforcing inventory limits in continuous time are worth the associated tractability hit.

\section*{Acknowledgments}
We acknowledge financial support from Air Liquide through the Center for Advanced Process Decision-making (CAPD) at Carnegie Mellon University. 
Akang Wang also gratefully acknowledges support from the James~C. Meade Graduate Fellowship and the H.~William and Ruth Hamilton Prengle Graduate Fellowship at Carnegie Mellon University, as well as the National Natural Science Foundation of China (Grant No. 12301416), the Shenzhen Science and Technology Program (Grant No. RCBS20221008093309021), and the Guangdong Basic and Applied Basic Research Foundation (Grant No. 2024A1515010306).

\bibliographystyle{plainnat}
\bibliography{bibliography}

\newpage
\appendix
\section{Nomenclature}
\label{sec:appendix}
\begin{tabular}{ p{2cm}   p{13cm} }
	\multicolumn{2}{l}{\bf{Parameters}} \\
         $c_{ij}^d$   &  cost to traverse the arc from node $(i, \alpha)$ to node $(j, \beta)$ by a vehicle in mode $d$ \\
         $T_{ij}^d$   &   time to traverse the arc from node $(i, \alpha)$ to node $(j, \beta)$ by a vehicle in mode $d$  \\
		$[W^{\ell}_{j\beta}, W^u_{j\beta}]$  & time window during which a vehicle may arrive at each customer $j \in V_c$ for the $\beta^\text{th}$ visit \\
         $H$   &   planning horizon  \\
         $K$          & homogeneous fleet size (number of vehicles) \\
         $Q$  &  vehicle capacity  \\
         $I_i^0/I_i^\ell/I_i^u$  &  initial/minimum/maximum inventory level at customer $i$ \\
         $r_i$   &  consumption rate for customer $i$  \\
		$\vartheta_i$ & minimum amount of product that the distributor has to ship to each customer $i \in V_c$ during the planning horizon\\
		$m_i$ & minimum number of visits at customer~$i \in V_c$ that have to be performed\\
	
	\multicolumn{2}{l}{\bf{Binary Variables}} \\
	 $x_{i\alpha j \beta} $   & equal to $1$, if the arc from $(i, \alpha)$ to $(j, \beta)$ is traversed in mode $d$; $0$ otherwise \\
        $y_{j \beta} $        & equal to $1$, if node $(j, \beta)$ is visited during the planning horizon; $0$ otherwise  \\

    \multicolumn{2}{l}{\bf{Continuous Variables}} \\
         $a_{j \beta}  $   &   time of arrival at node $(j, \beta)$, if visited \\
         $d_{j \beta}  $   &   time of departure from node $(j, \beta)$, if visited \\
         $\tilde{a}_{i \alpha j \beta}  $   &  time of arrival at node $(j, \beta)$, if visited immediately after node $(i, \alpha)$; $0$ otherwise \\
        $q_{j \beta}  $   &  total amount of product delivered to customer~$j$ at the $\beta^\text{th}$ visit \\     
		$f_{i \alpha j \beta}  $   &    amount of product transported from node $(i, \alpha)$ to node $(j, \beta)$ by a vehicle in mode $0$ \\  
         $\ell_{j \beta}  $   &    amount of product that a vehicle, after finishing its current trip, reloads at the depot before immediately proceeding to node $(j, \beta)$ \\
\end{tabular}

\newpage

\section{Validity of Rounded Capacity Inequalities}
    Given any subset of customers $S \subseteq V_c$, we can write
    \begin{equation}
        \begin{aligned}
            \sum_{j \in S}  \sum_{ \beta \in N_j } \left[ Q \sum_{(i, \alpha) \in \delta_0^+(j, \beta): \atop  i \notin S}    x^0_{i \alpha j \beta} +  Q \sum_{(i, \alpha) \in \delta_1^+(j, \beta) }  x^1_{i \alpha j \beta} \right]  \geq  \sum_{j \in S}  \vartheta_j,
        \end{aligned} \label{eq:RCI_demand}
    \end{equation}
    where the right-hand side represents the minimum amount of product required during the planning horizon, while the left-hand side represents the maximum amount of product delivered to customers in the set $S$; here, the second term accounts for the case when vehicles return back to the depot for replenishment and then immediately serve customers in the set $S$.
    
    Dividing both sides by $Q$ results in an RCI of the following form:
    \begin{equation}
        \begin{aligned}
            \sum_{j \in S}  \sum_{ \beta \in N_j } \left[  \sum_{(i, \alpha) \in \delta_0^+(j, \beta): \atop  i \notin S }    x^0_{i \alpha j \beta} + \sum_{(i, \alpha) \in \delta_1^+(j, \beta) }  x^1_{i \alpha j \beta} \right]  \geq  \ceil[\Bigg]{ \frac{\sum_{j \in S}  \vartheta_j }{Q} }.
        \end{aligned}
    \end{equation}

\end{document}